\documentclass[journal]{IEEEtran} 

\usepackage{cite}
\usepackage[]{graphicx}
\graphicspath{{images/}}

\usepackage{amsmath}
  {

	\interdisplaylinepenalty=2500
  }

\usepackage[utf8x]{inputenc} 
\usepackage{algorithm}
\usepackage[noend]{algpseudocode}
\usepackage{amsfonts}
\usepackage{amssymb}
\usepackage{booktabs}
\usepackage[caption=false,font=footnotesize]{subfig}
\usepackage{stfloats}
\usepackage{enumitem}
\usepackage{tikz}
\usetikzlibrary{shapes,arrows}

\usepackage{hyperref}[6.83]
\hypersetup{colorlinks,
			linkcolor=[rgb]{.61,0,0.3},
			citecolor=[rgb]{.14,.47,.14}}

\usepackage{tikz}
\usetikzlibrary{arrows.meta,
	calc, chains,
	quotes,
	positioning,
	shapes.geometric}

\newcommand{\norm}[1]{\left\Vert #1 \right\Vert}

\definecolor{OliveGreen}{rgb}{0.33, 0.42, 0.18}

\newcommand{\Hong}[1]{\textit{\color{blue}Hong: #1}}
\makeatletter
\newcommand{\removelatexerror}{\let\@latex@error\@gobble}
\makeatother

\definecolor{green2}{rgb}{0.0, 0.5, 0.0}

\newcommand{\review}[1]{{#1}}
\newcommand{\reviewtwo}[1]{{#1}}
\begin{document}

\title{Trust-Region Approximation of Extreme Trajectories in Power System 
Dynamics}
\author{D.~A.~Maldonado,~\IEEEmembership{Member,~IEEE,}
        E.~M.~Constantinescu, H.~Zhang, V.~Rao,
        and~M.~Anitescu,~\IEEEmembership{Member,~IEEE}
\thanks{This material is based upon work supported by the
	U.S. Department of Energy, Office of Science, Advanced Scientific
	Computing Research under Contract DE-AC02-06CH11357.}
\thanks{Authors are with the Mathematics and Computer Science Division at 
Argonne National Laboratory, Lemont, Illinois, U.S.A.}}%

\markboth{IEEE Transactions on Power Systems}%
{Optimization-Based Propagation of Uncertainties in Power System Dynamics}

\maketitle
\begin{abstract}
In this work we present a novel technique, based on a trust-region optimization 
algorithm and second-order trajectory sensitivities, to compute the extreme 
trajectories 
of power system dynamic simulations given \review{a bounded set that represents parametric uncertainty}. We 
show how 
this 
method, while 
remaining computationally efficient compared with sampling-based techniques, 
overcomes the limitations of previous 
sensitivity-based \reviewtwo{techniques to approximate the bounds of the trajectories} when the 
local 
approximation loses 
validity because  of the 
nonlinearity. We present 
several numerical experiments that showcase the 
accuracy and scalability of the technique,
including a demonstration on the IEEE New England test system.
\end{abstract}


\IEEEpeerreviewmaketitle 

\section{Introduction}
\label{section:intro}
\IEEEPARstart{T}he increase of complexity in  power distribution systems
is translating into uncertainty in transmission stability
studies \cite{Chaspierre2018}. Power distribution systems are traditionally 
modeled 
as lumped
loads at the transmission level; but as distribution becomes active,
the models become more complex \cite{NERC2016}, and their parameters become more 
difficult to
obtain from measurements or experiments. One must
understand how much the system dynamic response can vary given
the uncertainty of the system parameters. \review{
Three main types of models  capture parametric uncertainty in power system dynamics.

\begin{enumerate}
    \item When parameters vary with time with some known statistical regularity, they can be modeled as forcing stochastic processes \cite{Milano2013, Dhople2013, Adeen2021}. An example  is the load demand, whose short-term variation is observed to follow an Ornstein–-Uhlenbeck process \cite{Roberts2016}. This stochastic model turns the classical differential-algebraic equation (DAE) into a stochastic differential-algebraic equation (SDAE) that can be solved by using numerical techniques \cite{Dong2012}. 
    \item When parameters do not vary over time or their temporal variation is much slower than the transient response, we can represent the uncertainty by means of a probability distribution function \cite{Hockenberry2004, Xu2019}. Examples  are the parameters of a generator \cite{jovica}. In many cases, obtaining sufficiently descriptive data to create probability distributions is not feasible,  particularly when considering the parameterization of load models \cite{osti_1644285}.
    \item When parameters do not vary over time and their statistical properties are unknown, it is common to use a range or bounded set for the parameters  \cite{Hiskens2006, Kim2020}. This uncertainty model is also named \textit{unknown-but-bounded} \cite{Schweppe_Fred_C1973}. For this type of model, it is often of interest to compute the closure of all possible trajectories or extreme trajectories. Indeed, reliability criteria often specify allowable ranges for transient performance \cite{lesieutrereport, nerccriteria}.
\end{enumerate}}

\review{In this work we adopt the last type, and we represent uncertainty 
as a bounded set. We consider the problem of 
obtaining the extreme trajectories of the DAEs that model the power system in transient dynamics studies. \reviewtwo{While characterizing these extreme trajectories of the power system response can lead to conservative insights (i.e., the trajectory extremes might have a low probability of occurring), this type of analysis is useful for propagating parametric uncertainty in load models. Because load models are aggregated representations, obtaining statistical data of the parameters is challenging. Instead, parameter ranges are usually considered \cite{tenza2016analysis}.}

To solve this problem, one can use a Monte Carlo technique  that consists of sampling from the 
uncertainty set until statistical convergence \cite{jovica, Timko1983}. However, while Monte Carlo sampling} is general and robust since it makes no assumptions 
on the underlying model, its slow convergence makes it impractical for large 
parameter spaces. Other approaches to solve this problem are based 
on set-theoretic methods \cite{Zhang2020}, but these often have issues scaling 
to large systems. 

\review{More scalable approaches to characterize uncertainty have been developed using trajectory sensitivities.} 
As shown in \cite{Hiskens2006}, by constructing
local approximations of the DAE
solutions subject to perturbations of their parameters, one  can  
handle uncertainty efficiently. Inspired by the use of sensitivities in power 
systems,
more work has been done to improve the computational efficiency of the
trajectory sensitivities \cite{Zhang2017}. In particular, new
developments in the computation of second-order
sensitivities \cite{Geng2019} pave the way for efficient use of higher-order 
sensitivities in simulation codes.

The authors in \cite{Choi2017} propose a
methodology to \reviewtwo{approximate the extreme trajectories given a bounded uncertainty set}. The methodology consists of
constructing second-order Taylor approximations of the solutions of the
DAE and then solving an optimization problem to obtain the 
minimum and maximum 
at each time step. \reviewtwo{Formulating this trajectory bounding problem}
as an optimization problem allows \reviewtwo{it}
to be solved with a semidefinite programming algorithm. Despite the 
computational advantages 
of this method, however,  little 
attention is given to the
limitations of the Taylor approximation and the pitfalls of the optimization
procedure. Indeed, substituting the DAE solutions by a Taylor expansion can 
result in incorrect results due to the 
nonlinearity of the response. In these cases, the solution obtained will not 
correspond to a critical point of the actual DAE trajectory.

In our work we build on \cite{Choi2017} and develop a novel algorithm based on 
a trust-region optimization
methodology that ensures that the approximate surrogate model will be 
sufficiently accurate and  able to effectively approximate the 
extreme trajectories when the Taylor expansion deviates from the actual system 
response. 
\review{
The main contributions of this paper are the  following:
\begin{itemize}
  \item We propose a novel trust region-based approach to compute the extreme trajectories in dynamic simulations. By leveraging this nonlinear optimization method, we are able to increase the accuracy of the extreme trajectory approximation with respect to prior techniques.
  \item We carry out our analysis using dynamic models commonly used in commercial simulation software: a detailed synchronous generator, a governor, an exciter with saturation, and a composite load formed by passive and motor load.
  \item We show how to appropriately derive the initial conditions for the first- and second-order sensitivity systems, which are critical for quantifying the uncertainty of the initial conditions. This information is applied to the induction motor load, which has initial conditions that depend on the uncertain parameter.
  \item We derive a set of formulas in matrix form to integrate the first- and second-order sensitivities using the backward Euler method. 
\end{itemize}
}

In Section \ref{section:methodology} we formulate the extreme trajectory problem
mathematically, and we detail a computational solution consisting of a trust-region algorithm with a second-order expansion. In Section
\ref{section:casestudies} we  present experimental results for three
scenarios: a two-bus system with passive load  illustrating the basic
methodology, a two-bus system with dynamic load and control devices   showing the benefits of this method when nonlinearity of the response
increases, and a larger 39-bus system example showing the
scalability of this method. In Section \ref{section:conclusion} we 
present our conclusions.

\section{Methodology}
\label{section:methodology}
The power system dynamic equations can be abstracted into a general parameterized DAE system \cite{Milano2010}:
\begin{subequations}
\begin{align}
\dot{x} &= f(x, y, p, t) \, \\
0 &= g(x, y, p, t),
\end{align}
\label{eq:dae}
\end{subequations}
where $x \in \mathbf{R}^n$ is the differential state vector, $y \in
\mathbf{R}^m$ is the algebraic state vector, $p \in \mathbf{R}^p$ is a vector of
parameters, and $t$ is the time variable. The function $f()$ generally describes
time-domain elements of the system (such as generators and motors), and the 
function
$g()$ describes the network and power balance in the frequency domain. For
mathematical clarity, in this section we consider a concatenated vector $z=[x^T;y^T]^T$ and express these equations in a
succinct form using a mass matrix $M$:
\begin{align}
M \dot{z} &= h(z, p, t) \,,
\label{eq:ode}
\end{align}
where $M = (I, 0; 0, 0)$. The vector of parameters $p$ includes system 
parameters such as machine inertia and
load coefficients. In many cases, knowledge of $p$ is not exact, and it
might come  in the form of  either  a range or  a probability distribution
function. Often, one would ask the following question: \textit{What would be the
extreme trajectories of a certain dynamic quantity of interest, given that $p$
can take some predefined range of values?}

Given the nonlinear, discontinuous nature of the power system
model, this question has no straightforward answer. Mathematically,
\reviewtwo{finding the extreme upper trajectory} amounts to finding the argument of 

\begin{equation}
\label{eq:opt_problem}
\begin{aligned}
& \underset{p}{\text{maximize}}
& & {z}_i(p, t) \\
& \text{subject to}
& & p \in \Theta,
\end{aligned}
\end{equation}
for $t\in(t_0,t_N)$. \reviewtwo{Equivalently, the extreme lower trajectory is found by minimization}. Here, $\Theta = \{p \mid \theta_{\textit{min}} \leq p \leq 
\theta_{\textit{max}}\}$, and $z_i$ is a component of the state trajectory of 
the 
system. In general we cannot find closed solutions to the power 
system dynamic equations; hence, we need to either \textit{sample} or 
\textit{approximate}.

A common sampling approach is to use Monte Carlo techniques. Here, our model is
assumed to be a black box, and $\theta$ is sampled until statistical 
convergence. The benefits of this 
technique are that it can be implemented
nonintrusively on top of any existing simulation code, and it
also can be easily parallelizable. 
On the other hand, while much research on  efficient  sampling  methods  (importance  sampling,  quasi--Monte  Carlo,  etc.)  exists,  as the dimension of the parameter space increases, a greater number of samples are required to find an accurate solution, thereby rendering such methods intractable for large-scale systems.

Another approach is to approximate the equations over some local
region. This is accomplished with the use of derivatives, polynomial
interpolation, perturbation models, and  so on. For example, although an
analytical expression of $z_i(p, t)$ might not exist, we can
approximate it around a neighborhood $p_m$ with a Taylor
expansion: 
\review{
\begin{equation}
z_i(p_m + s, t) \approx z_i(p_m, t) + s^T \left . \frac{\partial z_i}{\partial p} 
\right |_{(p_m, t)}  + \frac{1}{2} s^T \left .\frac{\partial^2z_i}{\partial p^2} 
\right |_{(p_m, t)} s   \,.
\end{equation}
}

Here the derivatives \review{$\frac{\partial z_i}{\partial p}$ and $\frac{\partial^2z_i}{\partial p^2}$} can be 
obtained by computing the first- and second-order sensitivities with 
respect to 
$p$ of the DAE (\ref{eq:ode}). The first-order sensitivity can be 
obtained by solving the tangent linear sensitivity equations:
\begin{equation}
\review{ \frac{d}{dt} \left( \frac{\partial z_i}{\partial p} \right) = \frac{\partial  f}{\partial  z_i}\frac{\partial  z_i}{\partial  p} + \frac{\partial  f}{\partial  p} }.
\end{equation}

For more details on the derivation of the sensitivity equations for DAE 
systems 
 readers are referred to the foundational paper \cite{Hiskens2000} 
and the 
recent extension to second order \cite{Geng2019} and to \cite{Zhang2017} 
for 
details on the 
discrete adjoint method. In Appendix B we have 
included a summary of the continuous first- and second-order sensitivity 
equations and \review{derived a matrix expression to integrate the equations using the backward Euler technique}.

\begin{figure}[h!]
\centering
\includegraphics[width=0.4\textwidth]{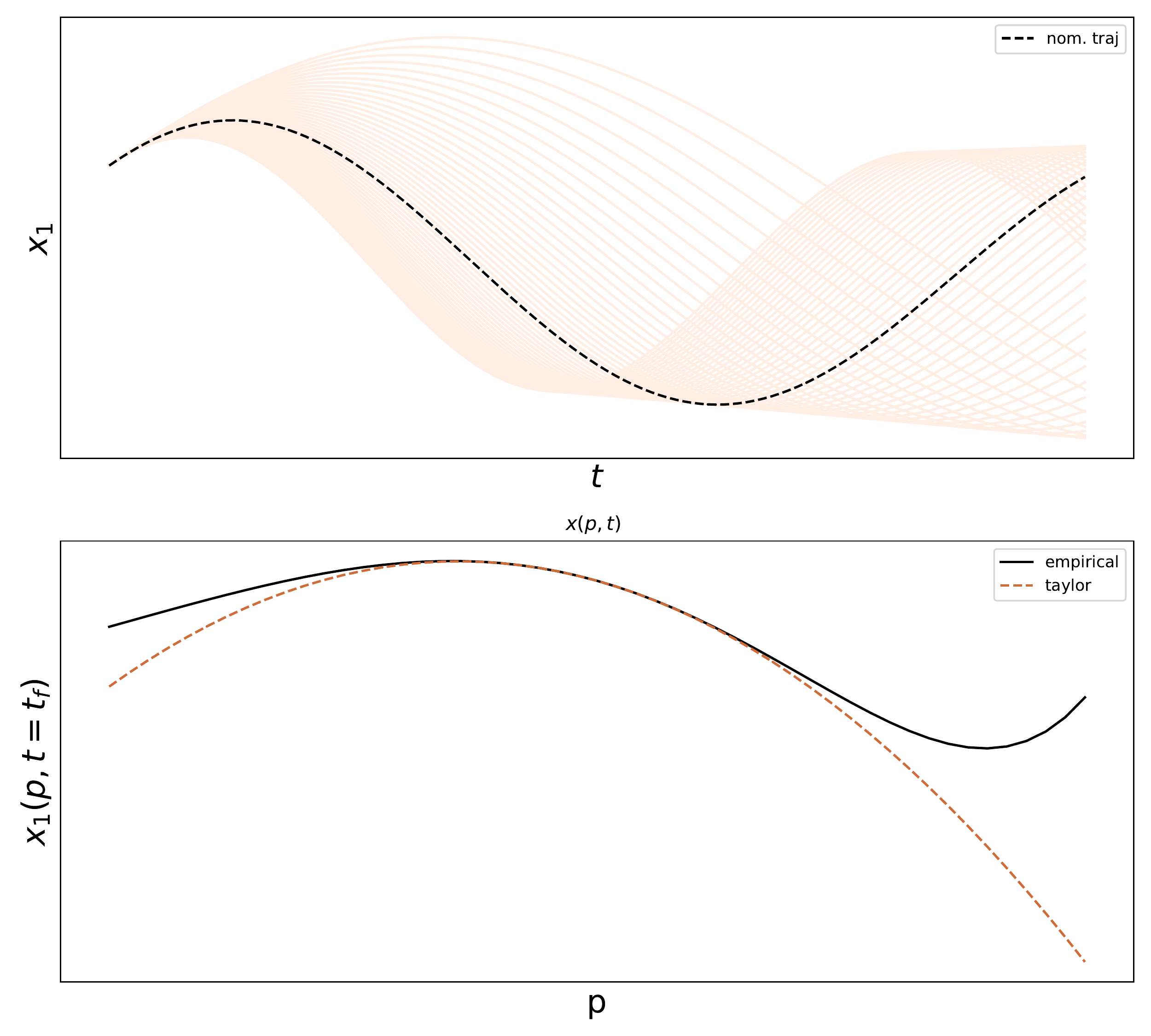}
\caption{On the top, multiple trajectories of a dynamical system given
variation of parameter $p$ are shown in orange, and the
black-discontinuous curve corresponds to the midpoint trajectory. On
the bottom, we \textit{slice} the trajectories on a certain time step
and plot the empirical function $x(t, p)$, to compare the empirical result (in black) with the second-order Taylor approximation (in orange).}
\label{fig:taylor_fail}
\end{figure}

Choi et al. \cite{Choi2017} use a local Taylor expansion as a \textit{surrogate
model} of $z(p, t)$ to solve (\ref{eq:opt_problem}). This is a good approach to
approximating the solution, but unfortunately the solution might be inaccurate
depending on the geometry of $z(p, t)$. For a certain range of the parameters, the
function may exhibit nonlinear characteristics that are not well approximated by
the quadratic form, and thus the calculated maximum and minimum
will not be
accurate. This situation  is illustrated in Figure \ref{fig:taylor_fail} showing how the Taylor
expansion at the midpoint does not capture the behavior of the function at the
minimum. 
Given that the trajectory functions (functions of the parameters) tend to be nonlinear and nonconvex, in our work we cast (\ref{eq:opt_problem}) as a nonlinear optimization problem, and we leverage the mature techniques in the numerical optimization field \cite{{Nocedal2006Numerical}}. Efficient techniques to solve these types of problems include line-search-based methods and trust-region-based methods. The availability of exact second-order information allows us to select methods with better convergence properties as opposed to methods that estimate the Hessian matrix of the cost function (e.g., quasi-Newton methods). In this work we chose trust-region methods to solve the optimization problem \eqref{eq:opt_problem}. Internally, trust-region methods use second-order Taylor approximations to the solutions of the DAEs. However, the key difference from the approach in \cite{Choi2017} is that the trust-region methods adapt the size of the trust-region based on how well the second-order Taylor series approximates the true DAE solution and will proceed iteratively until convergence. This approach aids in avoiding scenarios such as the one shown in Figure \ref{fig:taylor_fail} where optimizing the local Taylor expansion leads to wrong results. While the trust-region approach is an efficient and robust technique, we also note that similar performance could be achieved with line-search methods.

\subsection{Trust-region optimization}

 To solve an optimization problem of a general nonlinear function $f(x)$, 
 trust-region algorithms 
 \review{produce a
sequence of 
points $\{x_k\}$ that converge to a point $x_{*}$, where $\nabla f(x_*) = 0$ (i.e., $x_*$ is a first-order critical point)} by using 
a \review{quadratic} approximate 
surrogate model $m_k(x_k)$ of the objective function \cite{Conn2000}.
\review{This surrogate model is often a Taylor series expansion of $f(x)$ around $x_k$:
\begin{equation}
    m_k(x_k + s) = f(x_k) + \nabla f(x_k)^T s + \frac{1}{2} s^T \nabla^2 f(x_k) s \,,
\end{equation}
where $\nabla f(x_k)$ is the gradient of the objective function, $\nabla^2 f(x_k)$ is the Hessian, and $s$ is a small step around $x_k$.
The difference between $m_k(x_k + s)$ and $f(x_k + s)$ is small for small values of $||s||$. The optimization procedure consists of
reducing $f()$ through $m_k()$ by seeking a step $s$ that minimizes $m_k$. Depending on $s$, however, the approximation error can become
large enough that it will lead to incorrect results.}

The trust-region algorithm ensures that, at each iteration, the approximation is valid in a suitable neighborhood $\mathcal{B}_k$, where
\begin{equation}
	\mathcal{B}_k =  \{x \in \mathbf{R}^n \mid \norm{x - x_k} \leq 
	\Delta_k\}
	\,.
\end{equation}
This neighborhood is the \textit{trust region}. We compute a trial step $s_k$ to a trial point $x_k + s_k$ that 
minimizes the surrogate model $m_k()$ and satisfies $\norm{s_k}_2 \leq 
\Delta_k$. If the reduction of the surrogate 
model is in poor agreement with the reduction of the 
objective function $f()$, the trial point will be rejected, and the trust 
region will be reduced. We check the agreement with the following ratio:
\begin{align}\label{eq:trust_ratio}
\rho_k = \frac{\textit{actual reduction}}{\textit{predicted reduction}} = 
\frac{f_i(x_k) - f_i(x_k + s_k)}{m_k(x_k) - m_k(x_k + 
s_k)} \,.
\end{align}
The following algorithm \cite[Algorithm 4.1]{Nocedal2006Numerical} 
describes the basic trust-region method. Here, $\eta_1$ and $\eta_2$ are constants, and typically they are set to $\frac{1}{4}$ and $\frac{3}{4}$, respectively \cite{Nocedal2006Numerical}. We note that $\rho_k < 0$ corresponds to the case when the model predicts a decrease in function value at $m_k(x_k + 
s_k)$, but the function value,  $f_i(x_k + s_k)$, actually increased. In such scenarios, the step $s_k$ is rejected.
\begin{algorithm}
    \caption{Trust Region}
    \label{trust_region}
    \begin{algorithmic}[1] 
        \Procedure{Trust Region}{}
        \State Given $\hat{\Delta} > 0$, initialize $\Delta_0 \in (0, \hat{\Delta})$, and $\eta \in \left \lbrack 0, \eta_1 \right)$
        \For{$k = 0, 1, 2, \cdots $}
            \State obtain $s_k$ by reducing $m_k(x)$
            \State evaluate $\rho_k$ from \eqref{eq:trust_ratio}
            \If{$\rho_k < \eta_1$}
            \State $\Delta_{k+1} = \frac{1}{4} \Delta_k$
            \Else 
                \If{$\rho_k > \eta_2$ and $\|s_k\| = \Delta_k$}
                \State $\Delta_{k+1} = \textrm{min}(2 \Delta_k, \hat{\Delta})$
                \Else
                \State $\Delta_{k+1} = \frac{1}{4} \Delta_k$
                \EndIf
            \EndIf
            \If{$\rho_k > \eta$}
            \State$x_{k+1} = x_k + s_k$
            \Else 
            \State $x_{k+1} = x_k$    
            \EndIf
        \EndFor{}
        \EndProcedure
    \end{algorithmic}
\end{algorithm}

\review{If the step is accepted, we create a new surrogate model $m_{k+1}$ by evaluating the gradient and Hessian at $x_{k+1}$}.

Usually, the 
most 
expensive step in the trust-region algorithm is solving the subproblem (that is, reducing
\eqref{eq:surrogate_model}). For more details about the various strategies to 
solve the subproblem and convergence results, we refer the reader to 
\cite{Nocedal2006Numerical}. \reviewtwo{Note that the Hessian can be indefinite}. In this work we use a subproblem solver 
from 
\cite{lindon_roberts} \reviewtwo{that uses a variant of the Steihaug--Toint method \cite{Gould1999} with provisions for such cases.}

\subsection{Approximation of extreme trajectories}

To approximate the extreme trajectories, we use the trust-region 
approach to 
solve \eqref{eq:opt_problem}. In our case the objective function is a 
trajectory of the state variable at some time $t$ depending on a set of 
parameters $p$: $z^t_i(p)$, whereas the subproblem solved in 
the $k$th iteration can be written as follows (we omit the 
time superscript $t$ for brevity):

\begin{align} \label{eq:surrogate_model}
	\displaystyle \underset{d}{\textrm{min}} 
	\quad m_k(p_k + s) =  z_i + u_i^T s + \frac{1}{2}s^T 
	V_i s \,, \quad \textrm{s.t } \|s\| \leq \Delta_k\,, 
\end{align}
where $m_k(x)$ is a surrogate  based on a local quadratic 
approximation of the trajectory $z_i$ at $p$ and where $u_i$ and $V_i$ are 
respectively the first-order and second-order sensitivities of the 
$i$th state variable
evaluated at $p$, as described in Appendix \ref{section:trajimp}. 
Given the DAE system (\ref{eq:dae}) over a time interval $t \in (t_0, 
t_{\textit{end}})$, we discretize the interval into $j$ subintervals. Thus, 
the solution of (\ref{eq:opt_problem}) involves the solution of $j$ 
trust-region problems. For each time step $t_j$ the sketch of the computational 
procedure is as follows.

\begin{itemize}
\item Choose a point on the parameter interval: $p_\textit{nom} \in \Theta$. By 
default, we chose $p_\textit{nom} = p_m + \frac{p_M - p_m}{2}$.
\item Integrate the DAE \eqref{eq:dae} and its sensitivity equations from $t_0$ to $t_j$. With the results we can evaluate the surrogate model defined in (\ref{eq:surrogate_model}).
\item Proceed with Algorithm 1, and solve the subproblem. If it is necessary to 
evaluate another point $p_k$, this will involve another integration of 
the DAE and its sensitivities from $t_0$ to $t_j$. 
\end{itemize}

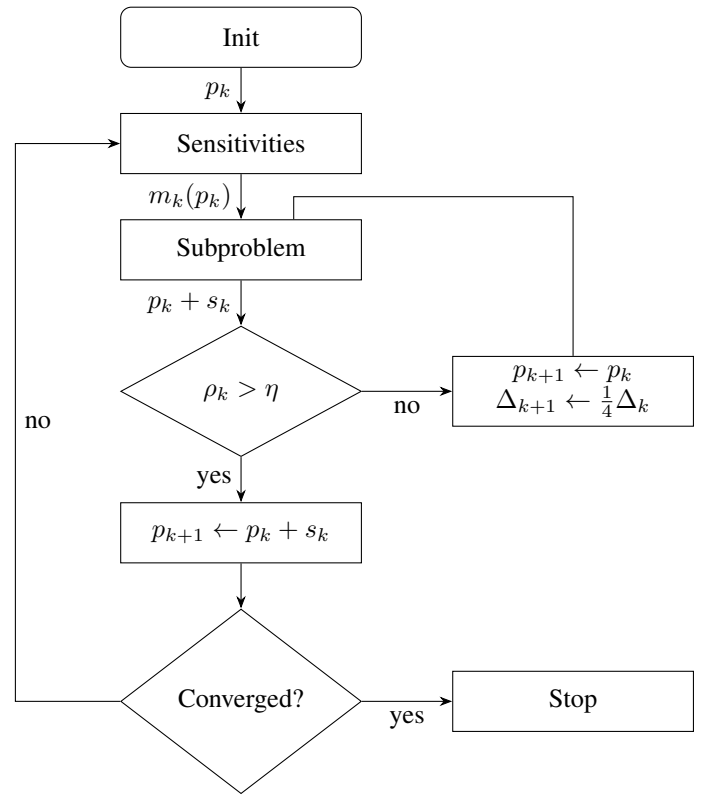
\begin{figure}
	\centering
 \begin{tikzpicture}[
	node distance = 6mm and 12mm,
	start chain = A going below,
	base/.style = {draw, minimum width=32mm, minimum height=8mm,
		align=center, on chain=A},
	startstop/.style = {base, rectangle, rounded corners},
	process/.style = {base, rectangle},
	io/.style = {base, trapezium, 
		trapezium left angle=70, trapezium right angle=110},
	decision/.style = {base, diamond},
	every edge quotes/.style = {auto=right}]
	]
	\node [startstop]       {Init};            
	\node [process]         {Sensitivities};
	\node [process]         {Subproblem};
	\node [decision]        {$\rho_k > \eta$};
	\node [process]         {$p_{k + 1} \leftarrow p_k + s_k$};
	\node [decision]        {Converged?};             
	\node [process,                             
	right=of A-4]    {$p_{k + 1} \leftarrow p_k$ \\ $\Delta_{k + 1} \leftarrow 
	\frac{1}{4}\Delta_k$};
	\node [process,                             
right=of A-6]    {Stop};
	\draw [arrows=-Stealth] 
	(A-1) edge["$p_k$"]          (A-2)
	(A-2) edge["$m_k(p_k)$"]    (A-3)
	(A-3) edge["$p_k + s_k$"]       (A-4)
	(A-4) edge["yes"]            (A-5)
	(A-5) edge[""]          (A-6)
	(A-4) edge["no"]          (A-7)       
	(A-7) |- ($(A-2.south east)!0.5!(A-3.north east)$)
	-| ([xshift=7mm] A-3.north)
	(A-6)  -- ++(-3,0) |- node[pos=.25, right]{no}   (A-2)
	(A-6) edge["yes"]     (A-8)
	;
\end{tikzpicture}
	\caption{Flowchart of the trust region. This represents the solution of a 
	problem for a particular time point $t$. We start $p_0 = p_{\textit{nom}}$, 
	and for each new $p_k$ we need to compute the sensitivities at the new 
	point.}
	\label{fig:diagram}
\end{figure}

This method will require repeated integration of the DAE system and can become 
more computationally expensive as $j$ increases. However, one  can  improve this situation. For instance, an initial integration of the DAE  and 
sensitivities with respect to $p_{\textit{nom}}$ can be done from 
$t_0$ to $t_{\textit{end}}$ and will serve as the initialization for 
the 
subproblem of the $k$ trust-region problems. Additional integrations will  
be carried out only if a trust-region problem requires an additional evaluation of the 
surrogate model. Furthermore, these extra computations can be performed 
in 
parallel. Figure~\ref{fig:diagram} shows the iterative process.

We have introduced this method as an \textit{approximation} of extreme 
trajectories. While our method would avoid situations like in 
Fig. 
\ref{fig:taylor_fail} where the second-order Taylor model does not effectively 
capture the geometry of the function over the parameter variation range, we 
could as well fall into a local minimum or maximum that would not correspond 
with an extreme trajectory point. We have tested our method with detailed 
dynamic models and with fault scenarios of wide magnitude ranges, and we have 
observed that the functions of the trajectories with respect to the 
parameters, 
while nonconvex,  do not seem to have multiple local minima.

In the next section we show experimentally 
that our method (a) improves the accuracy of \cite{Choi2017} when the system is 
stressed and exhibits acute nonlinear response and (b)  
approximates effectively the extreme trajectories, which we verify by comparing 
our results with sampling techniques.

\section{Case studies}
\label{section:casestudies}
One of the major sources of uncertainty in power systems is the
composition of the load. Since in normal scenarios the mixture of load
is not known with precision, we wish to understand how
system trajectories under a disturbance vary with the change of the
mixture combination \cite{Kim2020}. In other words, we would like to know how
sensitive the system is to a variation of the load mixture in a set of
buses. In most transient dynamic simulators, the load is modeled as
a current or power injection into the system---a Thevenin equivalent. For a 
certain bus the
active and reactive power injections can be written as
\begin{align}
P_{\textit{inj}}(t) &= P_1(t) + P_2(t) + \dots + P_n(t) \,, \\
Q_{\textit{inj}}(t) &= Q_1(t) + Q_2(t) + \dots + Q_n(t) \,,
\end{align}
where each $(P_i, Q_i)$ represents a distinct type of load. An example  
is the ZIP load model, where the individual power injections are passive 
functions of the voltage, and the 
composite load model, where the individual 
power injections have also dynamic load models such as the ones representing 
induction motor models and solar inverters. In this section we  obtain 
approximations to the extreme trajectories of generator states given 
perturbations of the 
components of the load mixture.

\subsection{Two-bus system with passive load}
\label{sec:casea}

As a first example we consider the load to be a mixture of constant impedance 
and constant power loads (a subset of the ZIP load model):
\begin{equation}
P^i_{\textit{inj}} = P_{\textit{z}} + P_{\textit{p}} \,.
\end{equation}
From the net load, $P^i_{\textit{inj}}$, a fraction, $\alpha$, will 
correspond to the constant impedance load and the 
rest of it, $1 - \alpha$, to the constant power load. A typical way to 
do this 
is by feeder usage statistics, which can be imprecise \cite{NERC2016}. To 
include this 
information in a dynamic simulation, we can write the constant impedance load as
\begin{equation}
P^i_{\textit{inj}} = \alpha \left(  \frac{V_i}{V_0} \right)^2 P_0 + (1 - 
\alpha) P_0  \,.
\end{equation}

\begin{figure}[h]
\centering
\includegraphics[width=0.45\textwidth]{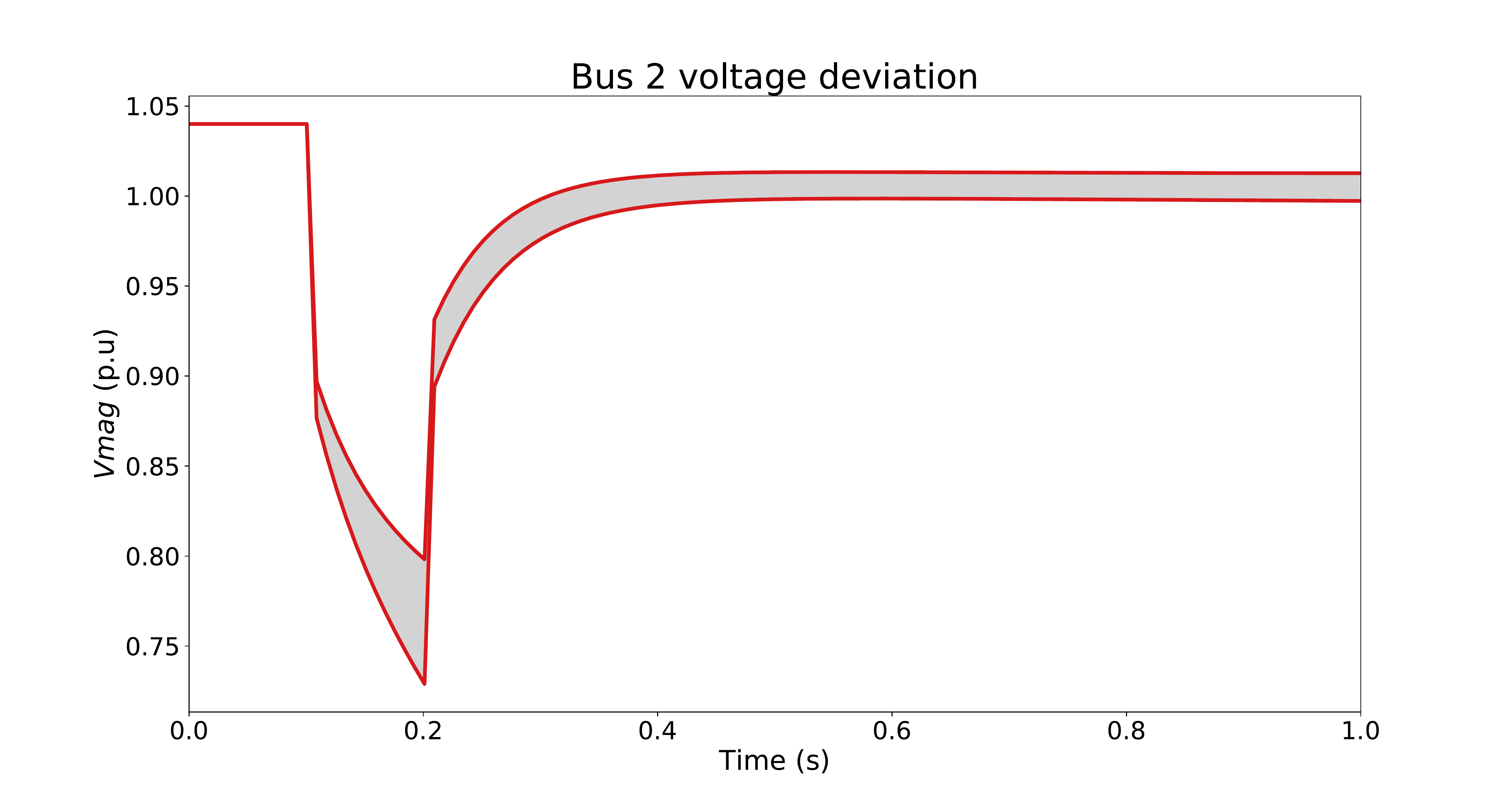}
\caption{In red, the extreme trajectories of the voltage magnitude at bus 2 
computed with the trust-region method. \review{In grey}, trajectories obtained from 
sampling the interval of $\alpha$ using a grid with $100$ points.}
\label{fig:casea_voltage}
\end{figure}

\begin{figure}[h]
\centering
\includegraphics[width=0.45\textwidth]{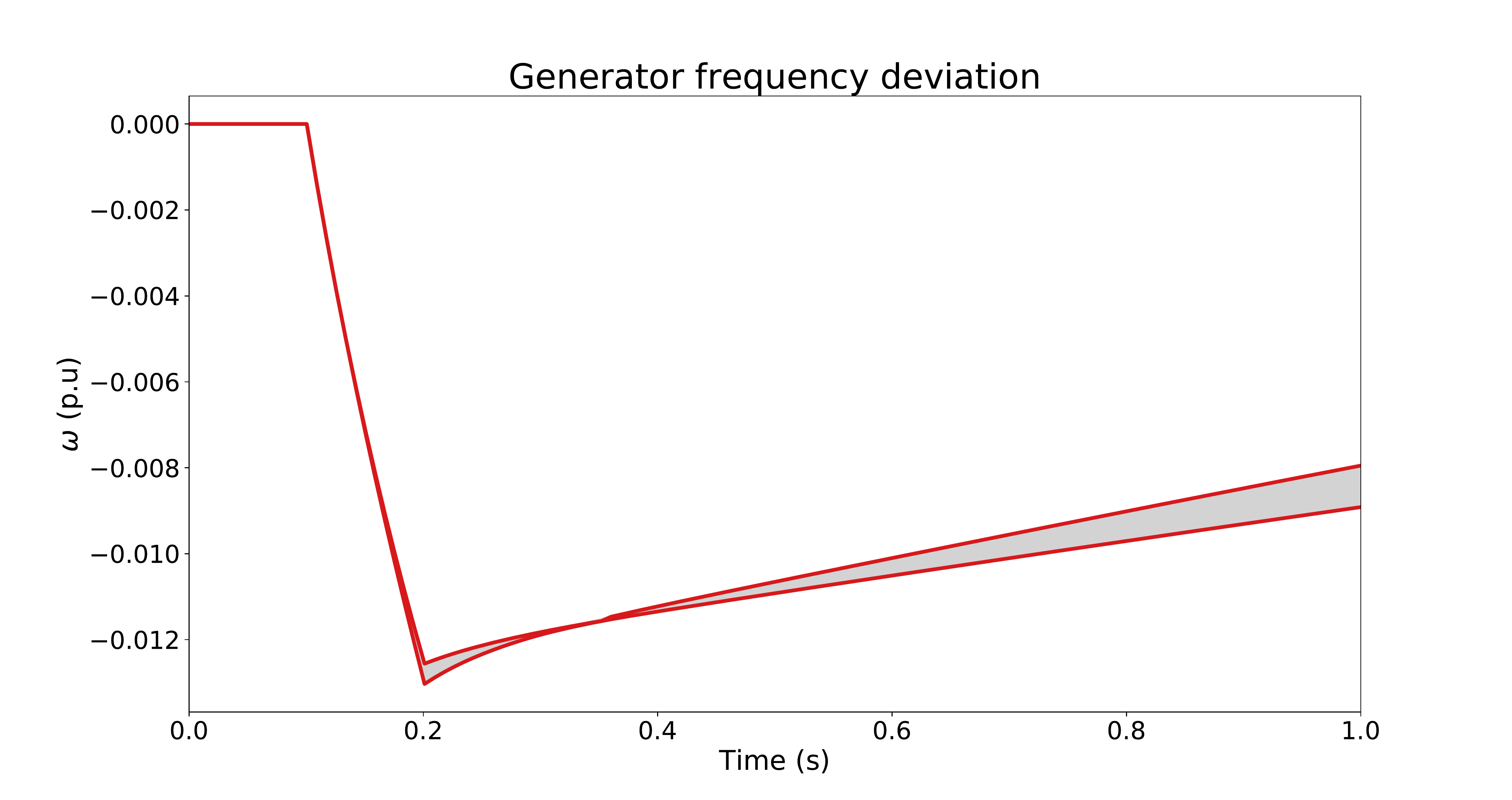}
\caption{In red, the extreme trajectories of the generator frequency deviation 
computed with the trust-region method. \review{In grey}, trajectories obtained from 
sampling the interval of $\alpha$ using a grid with $100$ points.}
\label{fig:casea_freq}
\end{figure}

\begin{figure}[h]
\centering
\includegraphics[width=0.45\textwidth]{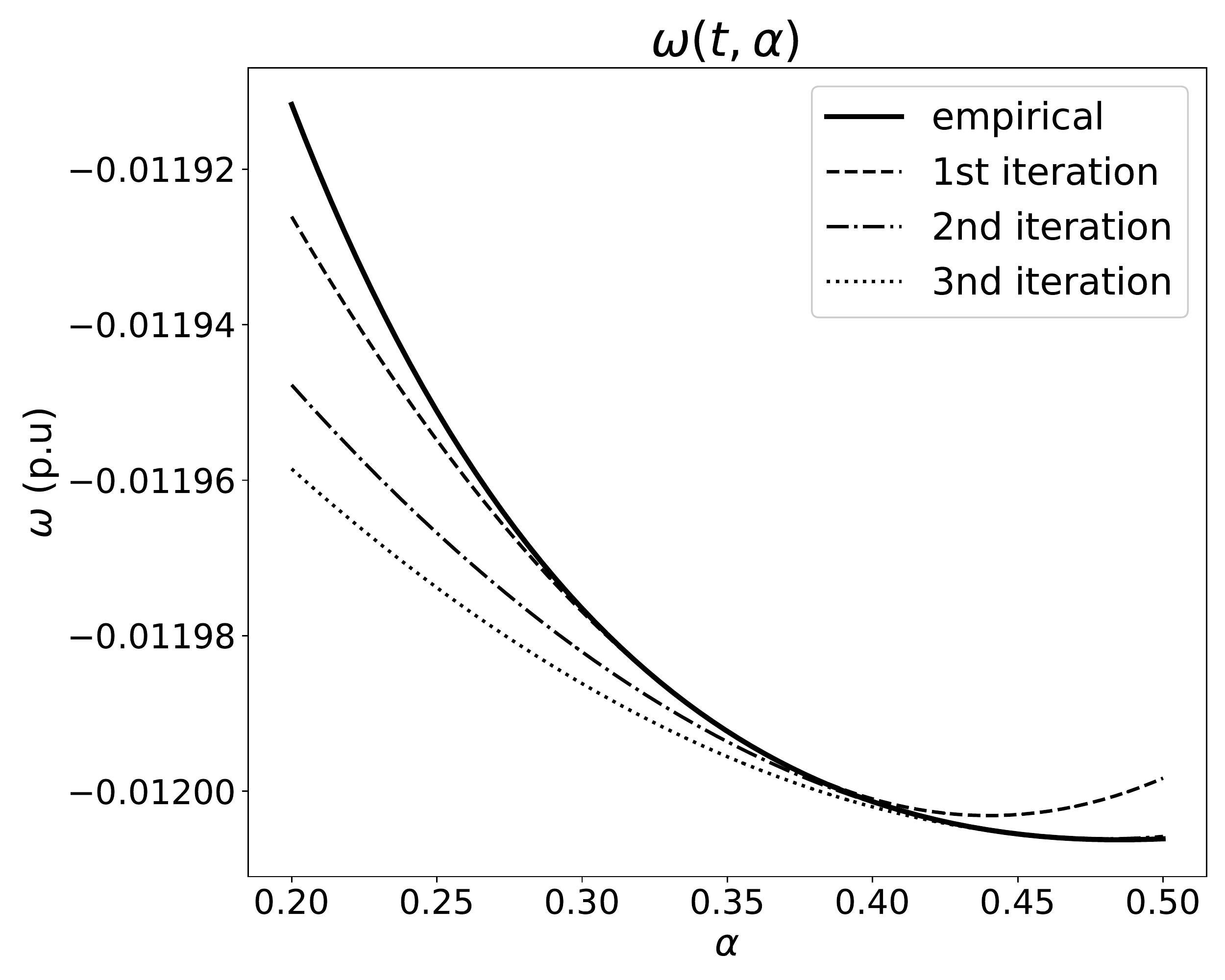}
\caption{The continuous curve is a plot of the empirical function $\omega(t, 
\alpha)$ at $t = 0.288$ s using Monte Carlo sampling. The \review{dashed} lines 
show 
the quadratic approximation at different iterations of the trust-region 
algorithm. We can see how the last iteration closely matches the minimum of the 
function.}
\label{fig:casea_taylor}
\end{figure}

We simulate a two-bus system consisting of a GENROU generator
model connected to this load. We 
would like to obtain 
the extreme trajectories of
the states of this system after a $0.6$ p.u. fault applied at $t =
0.1$ sec and removed at $t = 0.2$ sec, given that $\alpha \in [0.2,
  0.5]$. Figure~\ref{fig:casea_voltage} shows the computed extreme
trajectories for the voltage magnitude of the load bus. While the
effect of the parameter $\alpha$ on the variation of the voltage is
notable, its relationship is quasilinear and can be represented
effectively with the quadratic model \eqref{eq:surrogate_model}. This is not 
always the case, however; in
Fig. \ref{fig:casea_freq} we can see a similar plot of the generator
frequency deviation for a segment of time after the fault. While the
frequency variation is less conspicuous compared with the voltage variation, 
the frequency function
$\omega(t, \alpha)$ exhibits more nonlinearity, which 
results in
additional iterations. This can be seen
in Fig. \ref{fig:casea_taylor}, where the initial surrogate model
requires up to three iterations to approximate the function accurately.

By taking a \textit{slice} of the frequency trajectories ensemble at time $t = 
0.288$ s, we can observe how the trust-region algorithm works. To find the 
minimum, we first build a surrogate model around the nominal point $\alpha_1 = 
0.35$, and we solve the subproblem that finds a new trial step of the quadratic 
at 
$\alpha_2 = 0.364$. In this new point we compute again the function value and 
sensitivities and verify that $\rho > 0.75$, which allow us to accept the step 
and check that the gradient is still above tolerance. A new trial step with the 
new model at $\alpha_s$ gives us a trial step of $\alpha_3 = 0.4792$ where we 
find a critical point.  

We can define a simple metric to evaluate the accuracy of our method by 
computing the relative norm of the difference between the minimum computed with 
Monte Carlo sampling and the trust-region method:
\review{
\begin{subequations}
\begin{align}
	\epsilon_M &= \frac{\norm{M_\textit{mc} - M_\textit{trust}}}{\norm{M_\textit{mc}}} \,, \\
	\epsilon_m &= \frac{\norm{m_\textit{mc} - m_\textit{trust}}}{\norm{m_\textit{mc}}} \,.
\end{align}
\end{subequations}
}
The results for this example are given in \ref{table:casea}, where we find 
good agreement with the sampling-based method.

\begin{table}

	\label{table:casea}
	\renewcommand{\arraystretch}{1.3}
	\caption{Case A: approximation relative errors}
	\label{table_example}
	\centering
	\review{
	\begin{tabular}{c c c}
		\hline
		\bfseries Variable & \bfseries $\epsilon_M$ & \bfseries $\epsilon_m$ \\
		\hline
		$V_{\textit{mag}}$ & 5.736e-09 & 1.320e-09 \\
		$\omega$ & 5.049e-04 & 4.171e-07\\
		\hline
	\end{tabular}}
\end{table}

\subsection{Two-bus system with passive and motor load}
\label{sec:caseb}

We can modify the previous example to include an induction motor model instead:
\begin{align}
P_{\textit{inj}}(V_0, t_0) &=\alpha P_{\textit{z}} + (1 - \alpha) P_{\textit{mot}} \,, \\
Q_{\textit{inj}}(V_0, t_0) &= \alpha Q_{\textit{z}} + (1 - \alpha) Q_{\textit{mot}} \,.
\end{align}
The motor equations are
\begin{subequations}
\begin{align}
\dot{e'_{d}} &= -\frac{1}{T_p}(e'_d + (x_0 - x')i_q) + s \omega_s e'_q \,, \\
\dot{e'_{q}} &= -\frac{1}{T_p}(e'_q - (x_0 - x')i_d) - s \omega_s e'_d \,, \\
\dot{s} &= \frac{1}{2H}(\tau_m - e'_d i_d - e'_q i_q) \,, \\
0 &= r_a i_d - x' i_q + e'_d + V sin(\theta) \,, \\
0 &= r_a i_q - x' i_d + e'_q  - V cos(\theta) \,,
\end{align}
\end{subequations}
where $x = \begin{bmatrix} e'_{d} & e'_{q} & s \end{bmatrix}^T$ and $y = 
\begin{bmatrix} i_d & i_q \end{bmatrix}^T$ are the state and algebraic variables for the DAE system, respectively. 
The torque value $\tau_m$ is to be found. The active  power and reactive power  
consumed by the motor are written as
\begin{subequations}\begin{align}
P_{\textit{mot}} &= - V sin(\theta) i_d + V cos(\theta) i_q \,, \\
Q_{\textit{mot}} &=   V cos(\theta) i_d + V sin(\theta) i_q \,.
\end{align}\end{subequations}
We initialize the motor with $V_0, \theta_0$ and $(1 - \alpha) 
P_{\textit{inj}}$ from a power flow solution and obtain a set of initial 
states $z_0 = \begin{bmatrix} e'_{d} & e'_{q} & s & i_d & i_q & 
\tau_m\end{bmatrix}^T$. After setting the initial states (and 
parameter $\tau_m$), the motor reactive power consumption, in general, does not 
match $(1 - \alpha Q_0)$. To fix this discrepancy,  we introduce a shunt reactance:
\begin{equation}
y_{\textit{sh}} = \frac{V_0 cos(\theta_0) i_d + V_0 sin(\theta_0) i_q - \alpha Q_0}{(V_0)^2}.
\end{equation}

Computing the sensitivities with respect to $\alpha$ in the case where a motor 
is one of the components of the load presents additional complications. In this 
case, the initial states and parameters $\tau, y_{\textit{sh}}$ from the motor 
are all dependent on $\alpha$. Hence, we need to obtain the sensitivities of 
the initial state with respect to $\alpha$
$\frac{d\dot{x}}{d \alpha}$. Writing the system equations in the  form 
\eqref{eq:ode}, we can derive the sensitivity equations with implicit 
differentiation:

\begin{equation}
u^{\alpha} = -\mathcal{J}^{-1} H_{\alpha} \,.
\end{equation}
A similar derivation for the second-order sensitivities gives
\begin{equation}
v^{\alpha \alpha} = -\mathcal{J}^{-1} \left((I_m \otimes (u^{\alpha})^T) 
\mathcal{H} 
u^{\alpha} + 
2H_{x\alpha} + H_{\alpha \alpha} \right) \,,
\end{equation}
where $\mathcal{J}$ and $\mathcal{H}$ are the Jacobian and Hessian matrices of 
$\eqref{eq:ode}$,  respectively. More notation details are presented in Appendixes 
\ref{section:notation} and \ref{section:trajimp}. For a 
practical 
implementation, we have found it  more convenient to specify the static 
parameters $\tau, y_{\textit{sh}}$ as additional state variables and augment the original DAE with
\begin{align*}
\dot{\tau} &= 0 \,,\\
\dot{y}_{\textit{sh}} &= 0 \,.
\end{align*}

In addition to the induction motor, we  include a governor and 
an exciter with saturation. Their inclusion does not incur additional 
technical difficulties and allows us to obtain a system response richer in 
nonlinear features to better justify the effectiveness of the trust-region 
method. A \review{$0.05$} p.u. fault is applied at bus 2 from $t=0.25$ to $t = 0.40$ 
seconds; $\alpha$ is restricted in the range $[0.3, 1.0]$.

\begin{figure}[h]
	\centering
	\includegraphics[width=0.45\textwidth]{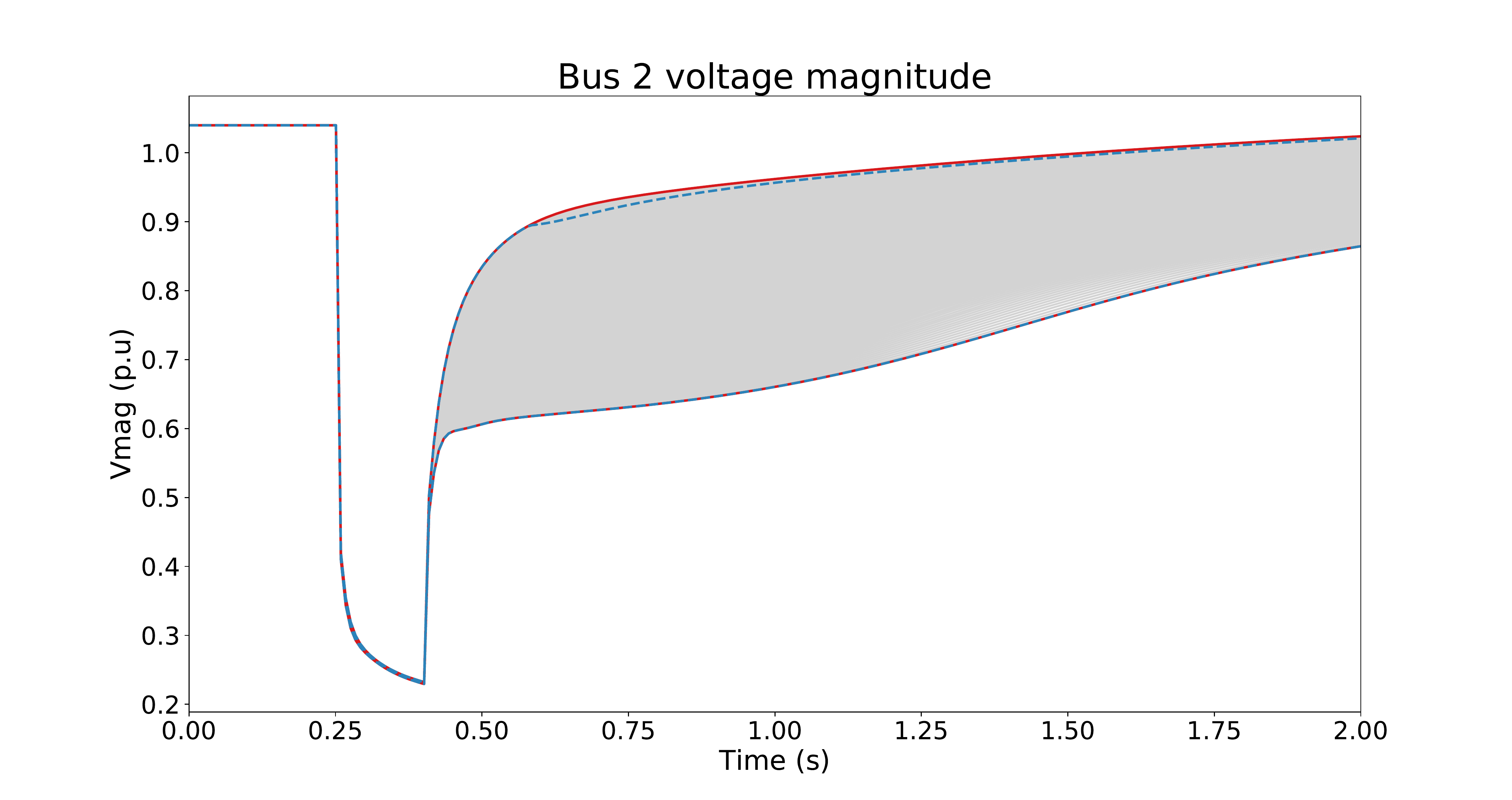}
	\caption{In red, the extreme trajectories of the voltage magnitude at bus 2 
		computed with the trust-region method. \review{In grey}, trajectories obtained 
		from sampling the interval of $\alpha$ using a grid with $1,000$ points. 
		In \review{dashed-}blue, the optimization is carried out with the 
		model at the nominal point.}
	\label{fig:caseb_traj_vmag}
\end{figure}

\begin{figure}[h]
	\centering
	\includegraphics[width=0.45\textwidth]{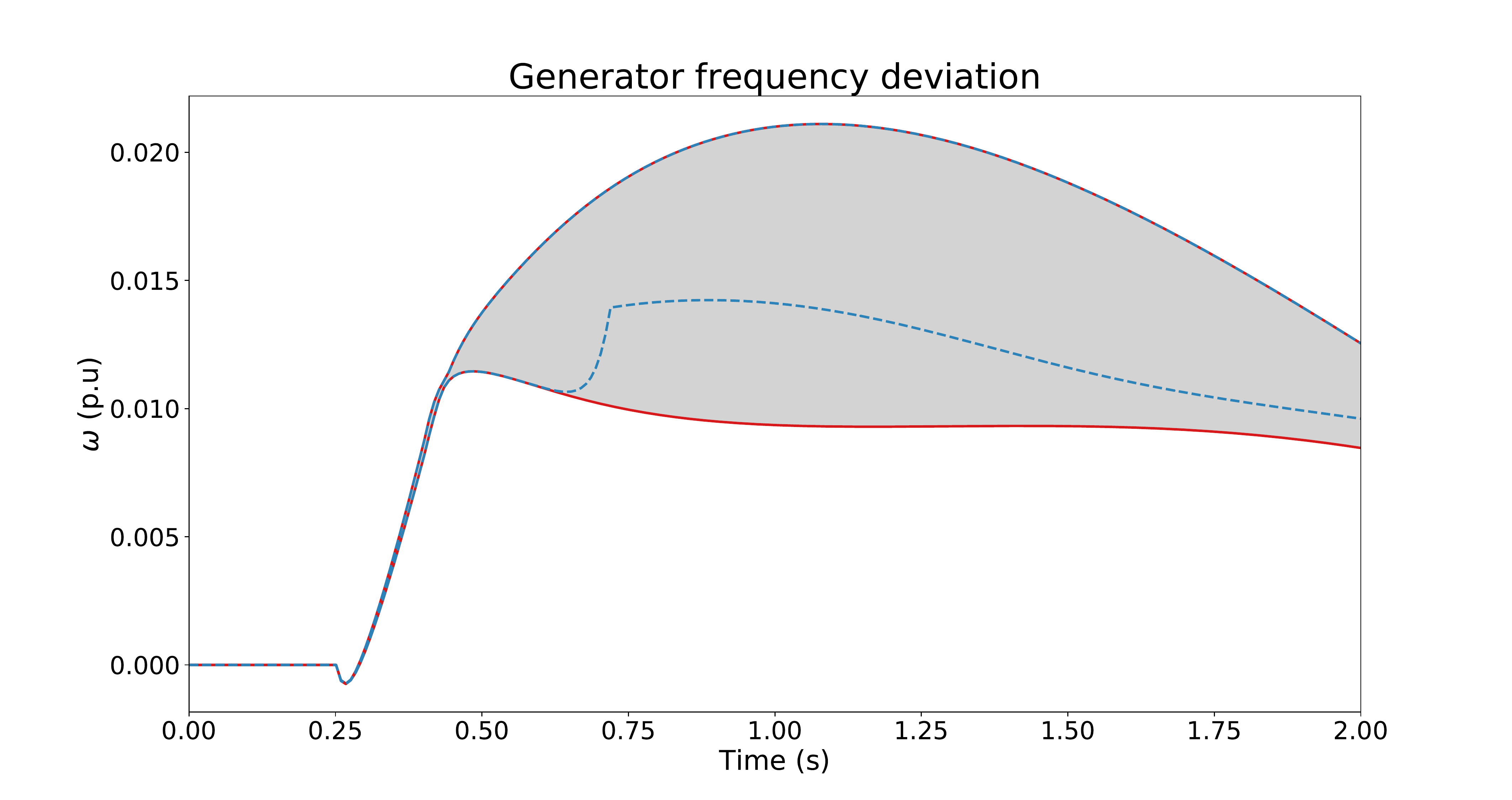}
	\caption{In red, the extreme trajectories of the generator frequency deviation 
		computed with the trust-region method. \review{In grey}, trajectories obtained 
		from sampling the interval of $\alpha$ using a grid with $1,000$ points. 
		In \review{dashed-}blue, the optimization is carried out with the 
		model at the nominal point.}
	\label{fig:caseb_traj_freq}
\end{figure}

\begin{figure}[h]
	\centering
	\includegraphics[width=0.45\textwidth]{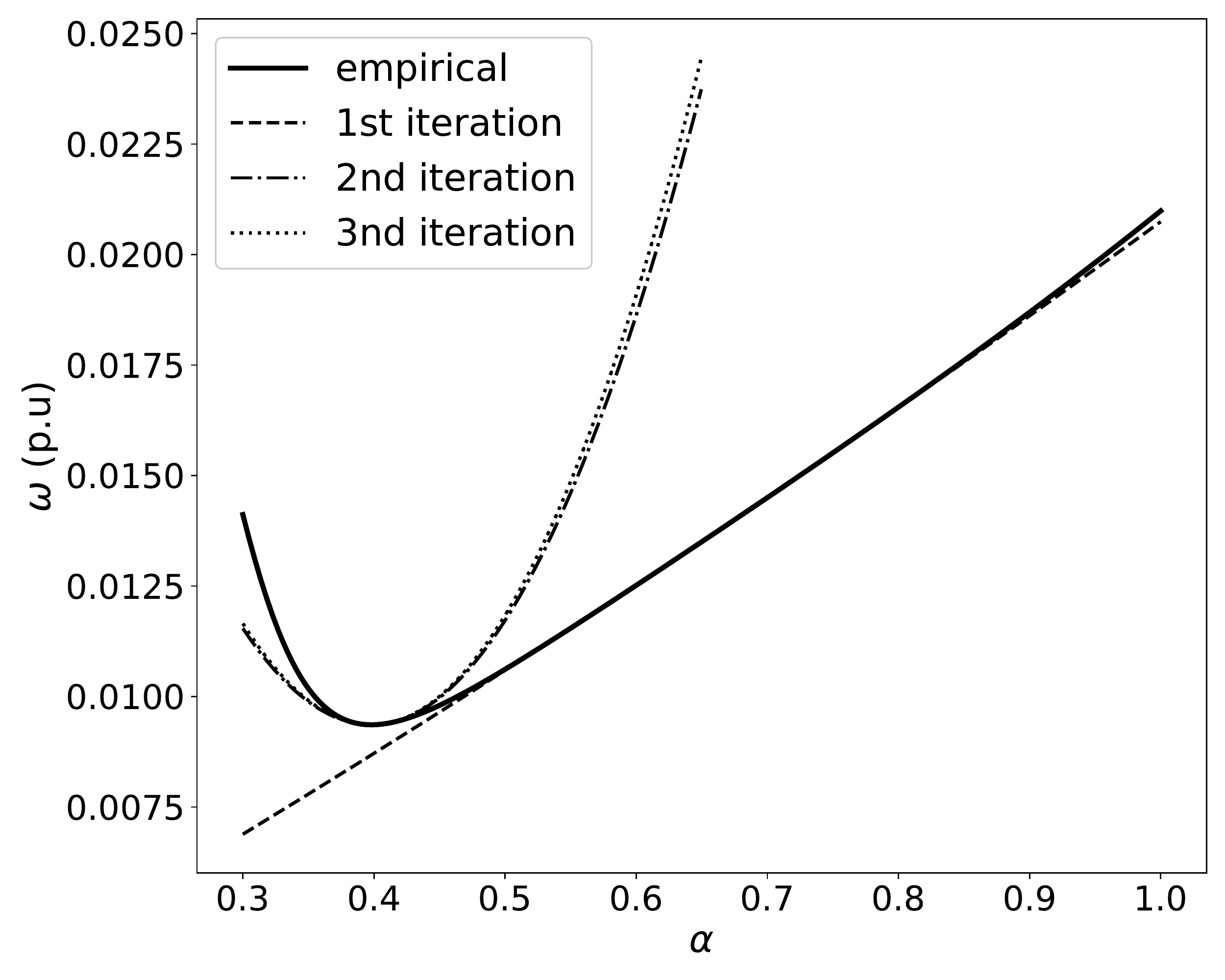}
	\caption{The continuous curve is a plot of the empirical function $\omega(t, 
		\alpha)$ at $t = 0.288$ using Monte Carlo sampling. The \review{dashed} lines 
		show 
		the quadratic approximation at several iterations of the trust-region 
		algorithm.}
	\label{fig:caseb_slice_freq}
\end{figure}

We compare the trust-region method with the approach delineated in
\cite{Choi2017} where only the Taylor expansion at the nominal point is used.
Figures~ \ref{fig:caseb_traj_vmag} and \ref{fig:caseb_traj_freq} show that the
trust-region approach can compute the extreme trajectories with higher accuracy,
particularly with regard to frequency response. The trust-region minimization
proceeds as follows: the surrogate model is built at the nominal point $\alpha =
0.65$, resulting in the trial point $\alpha = 0.3$; $\rho$ is found to be below
$\frac{1}{4}$, and thus the trust radius is reduced, which leads to a new trial
point at $\alpha = 0.4$. This iterative process, as illustrated in
Fig.~\ref{fig:caseb_slice_freq}, eventually converges to the critical point at
$\alpha = 0.398$. Figure~\ref{fig:caseb_slice_freq} also shows that the initial Taylor
model does not capture the curvature of the model. Quantified evaluation of the
errors for the two approaches is given in Table \ref{table_exampleb}.

\begin{figure}[h]
	\centering
	\includegraphics[width=0.45\textwidth]{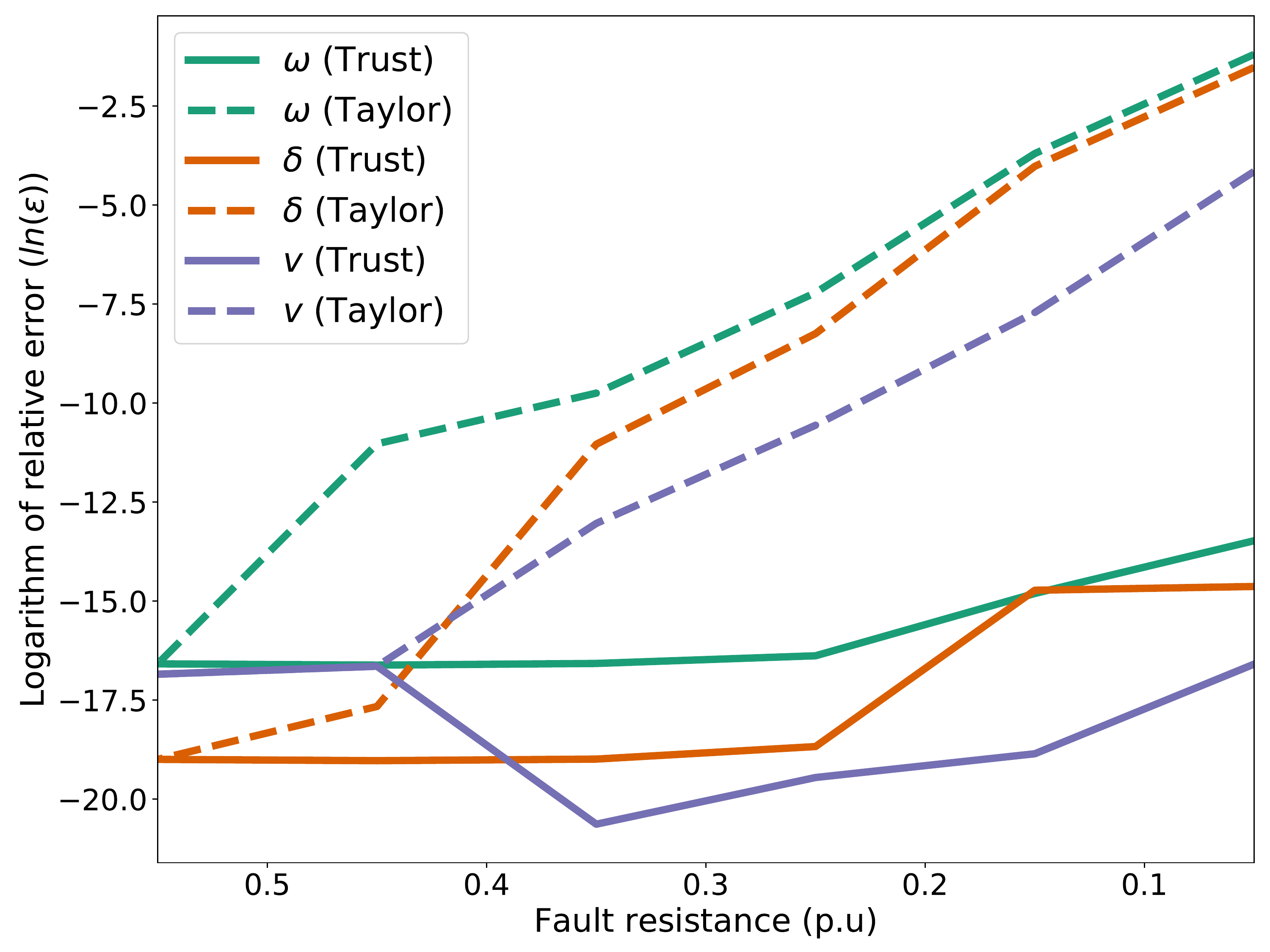}
	\caption{\review{Logarithm of the relative error for the Taylor (dashed line) and the trust-region (continuous line) algorithms for decreasing fault resistance values. The state variables depicted are the generator speed, the generator angle, and voltage magnitude. As we decrease the resistance of the fault, the perturbation of the system is greater, leading to large nonlinearities. We can see how the error incurred by our method remains controlled.}}
	\label{fig:caseb_log}
\end{figure}

\review{We next perform additional simulations with this same system, where we compute the relative errors
for a range of fault values from $0.55$ p.u to $0.05$ p.u. As the fault resistance decreases, the disturbance
is greater, which leads to a higher degree of nonlinearity. In Fig.~\ref{fig:caseb_log} 
we have plotted the error
of both the Taylor approach and the trust-region approach for a set of state variables. We can see that while these methods
incur  similar error for mild disturbances, as nonlinearity grows, the error of the trust-region method remains small.}

\begin{table}
	\renewcommand{\arraystretch}{1.3}
	\caption{Case B: approximation relative errors}
	\label{table_exampleb}
	\centering
	\review{
	\begin{tabular}{c c c c c}
		\hline
		\bfseries Variable & $\epsilon^{\textit{Taylor}}_M$ & $\epsilon^{\textit{Trust}}_M$
		& $\epsilon^{\textit{Taylor}}_m$ & $\epsilon^{\textit{Trust}}_m$\\
		\hline
		$V_{\textit{mag}}$ & 4.879e-03 & 4.172e-09 & 9.116e-07 & 8.407e-09\\
		$\omega$ & 3.991e-08 & 3.991e-08 & 0.303 & 1.406e-06\\
		\hline
	\end{tabular}}
\end{table}

\subsection{New England test system}

\begin{figure}[h]
	\centering
	\includegraphics[width=0.45\textwidth]{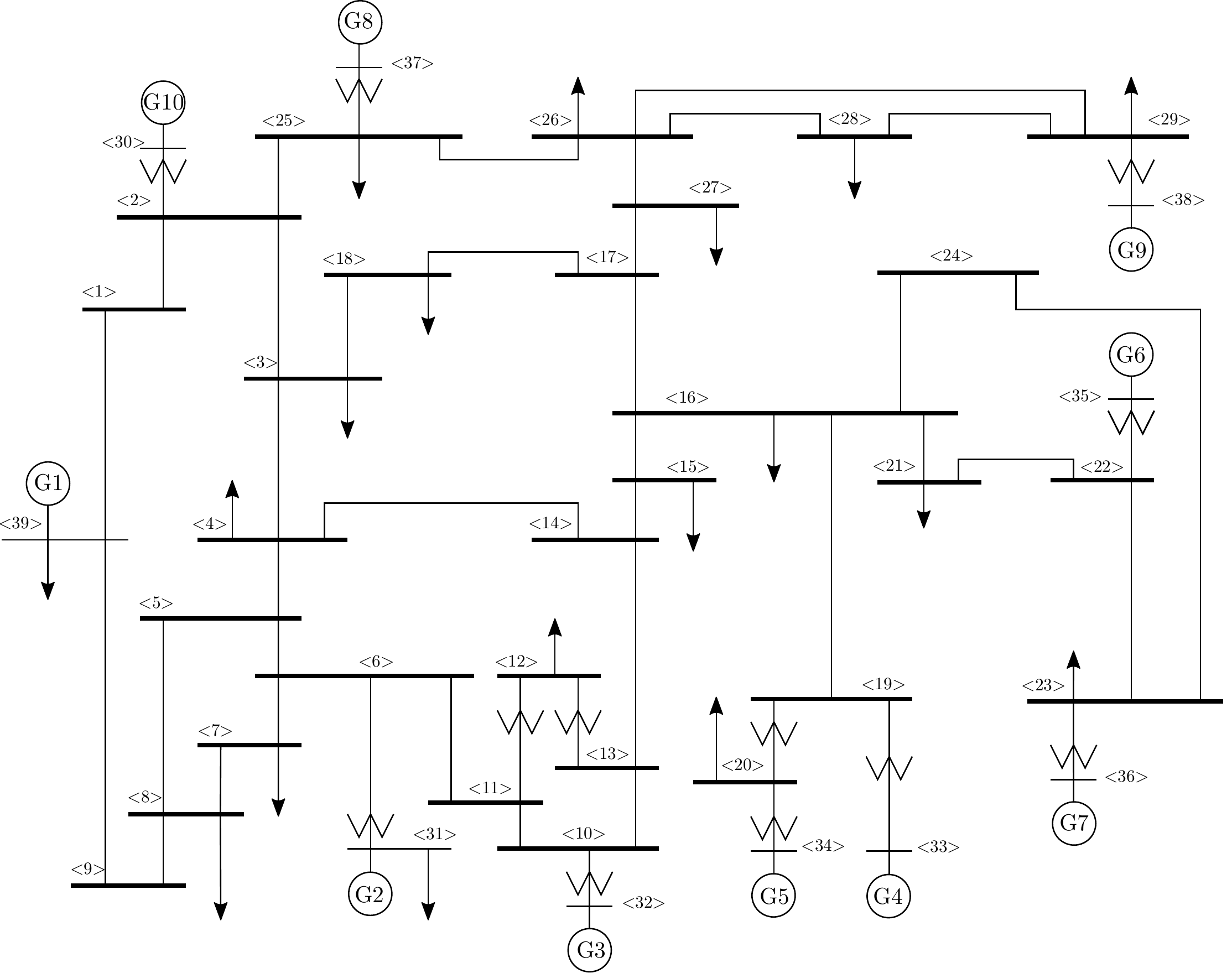}
	\caption{New England test system. A three-phase to ground fault is applied 
	at bus 2. We examine the behavior of generators 1 and 10 at buses 30 and 39, 
	respectively.}
	\label{fig:neng}
\end{figure}

In this last example we show the scalability of our algorithm with the New 
England test system (Fig.~\ref{fig:neng}). Here we revert to the passive ZIP 
load 
model 
introduced in 
\ref{sec:casea}. In the previous cases we had a single parameter, but now the 
dimension of the parameter space is 19 (that is, the number of loads), and 
computations become more expensive.

\begin{figure}[h]
	\centering
	\includegraphics[width=0.45\textwidth]{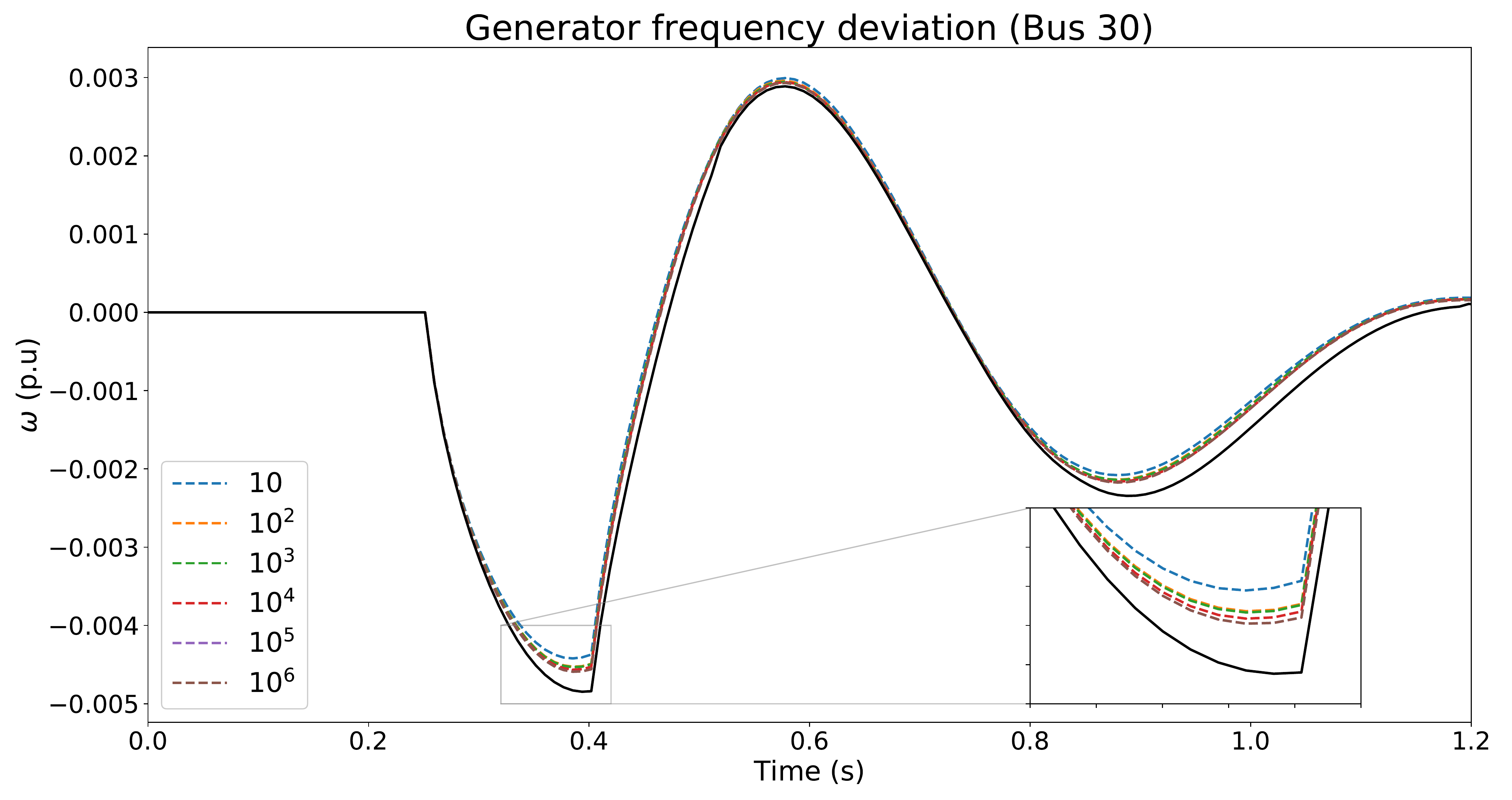}
	\caption{In black, the minimum trajectory computed with the trust-region 
	method. \review{Dashed} lines represent the same minimum trajectory computed 
	with Monte Carlo with increasing number of samples. We can see that Monte 
	Carlo converges slowly.}
	\label{fig:casec_gen30f}
\end{figure}

Because of the high dimensionality of the parameter space, we cannot 
plot the empirical function as we did in the preceding sections.
In Fig.~\ref{fig:casec_gen30f} we show the minimum trajectory of the frequency 
deviation of the generator at bus 30 given that for each load $k$ the parameter 
$\alpha_i \in [0, 1]$ using our trust-region technique. We also use Monte Carlo sampling
to obtain the minimum trajectory by sampling uniformly from the parameter space, 
and we plot the result for increasing amounts of samples.

\review{We can see in the} augmented 
section plot of Fig.~\ref{fig:casec_gen30f} how Monte Carlo (MC) converges slowly. We observe that, despite 
increasing the sampling number by an order of magnitude, subsequent MC 
simulations approach the trust-region solution with smaller increments each 
time. 

\begin{figure}[h]
	\centering
	\includegraphics[width=0.45\textwidth]{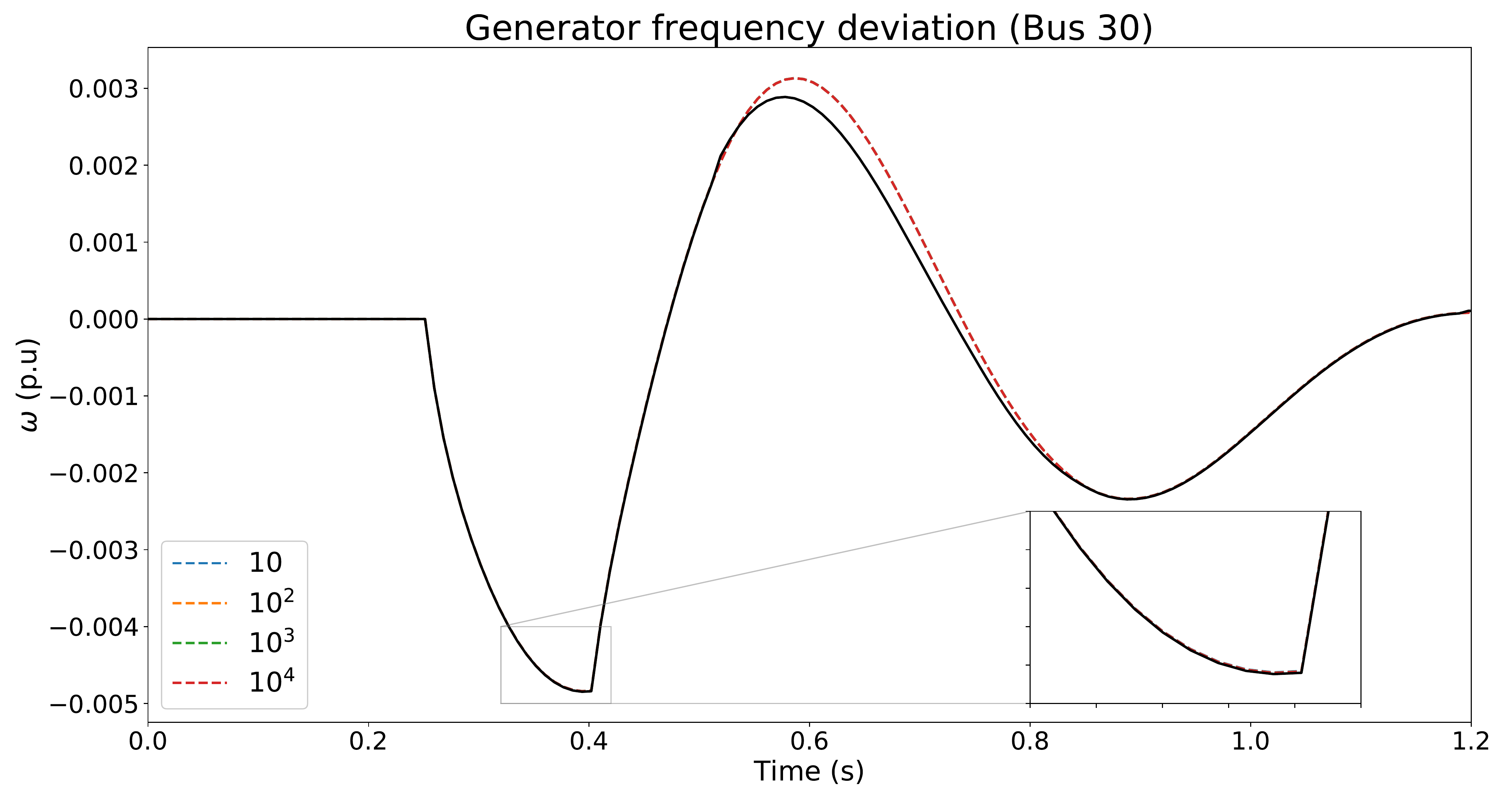}
	\caption{In black, the minimum trajectory computed with the trust-region 
		method. \review{Dashed} lines represent the same minimum trajectory 
		computed 
		with Monte Carlo with increasing number of samples. This time the 
		samples are drawn from a smaller region that we estimate to be closer 
		to the minimum in the interval $t = 0.25$ to $t = 0.42$.}
	\label{fig:casec_gen30f_s}
\end{figure}

To verify that indeed the trust-region solution that we show in 
Fig.~\ref{fig:caseb_traj_freq} is the minimum we would obtain with a 
sampling-based procedure, we sample around a small region
near the trust-region solution. For the nadir of Fig.~\ref{fig:casec_gen30f},
roughly from $t = 0.25$ to $t = 0.42$, we observe that the trust-region
algorithm finds a local minimum at one of the corners of the parameter space in
which $\alpha_i = 0$. Now we sample uniformly from the region $\alpha_i \in [0,
1]$ and observe that the trust-region solution closely matches the MC one (see
Fig.~\ref{fig:casec_gen30f_s}). \reviewtwo{Using the trust-region technique, we were able to compute this case in 2 minutes and 37 seconds using a prototype implementation in the Python programming language, whereas the computational time of MC is 17 minutes and 16 seconds for $10,000$ samples.}

In some trajectories we observe that the extreme boundaries of the parameter 
region (the 
corners) correspond to extremes of the function, and this situation would motivate the use 
of set-theoretic techniques. However, this is not true in 
general---especially when acute nonlinearities and discontinuities, typical of 
power system 
dynamics under large disturbances, arise in the simulation.

\section{Conclusion and future work} \label{section:conclusion}

In this work we have proposed an algorithm to compute trajectory extremes in 
power system dynamics using second-order sensitivities and a trust-region 
optimization approach. In addition to previous work where Taylor expansions 
were used as a surrogate of the DAE solution, our method guarantees that the 
local minimum computed with the sensitivities will be acceptably close to the 
actual DAE solution without a large increase in computational expense.

We have exemplified our technique with generator and load models used in common 
transient dynamics software packages (e.g., round rotor dynamic model with 
exciter, governor, and saturation) showing that the technique can be applied to 
realistic power systems. Computing the 
effect of load parameters on the dynamic trajectories is an important task as 
load models become more complex. For this purpose we have introduced a 
mixture model of a passive load and an induction motor and derived the first- 
and second-order sensitivities of the initial conditions of the motor with 
respect to the mixture parameter. 

We have also derived a backward Euler implementation of the first- and 
second-order
sensitivities computation, and we have detailed ways to compute these 
effectively.

In the future it will be of interest to 
extend this optimization method to models that present protection-induced 
discontinuities to analyze  phenomena such as fault-induced delayed 
voltage recovery. Furthermore, describing the power 
system with models of prespecified nonlinearity, such as the Lu\'re system 
formalism, 
could 
help us 
obtain closed-form results of higher generalization power. \reviewtwo{Our method has considerable room to be optimized such that it can be used for much larger systems. In the future we will work on efficient computation and usage of second-order sensitivities that include parallelization and randomization.}

\appendix

\subsection{Differentials and Notation} 
\label{section:notation}
Let $\mathbf{f}:S \to \mathbf{R}^m, S \in \mathbf{R}^n$ be a vector 
function, 
and let z be an interior point of $S$. If we can find a matrix 
$\mathcal{J} \in 
\mathbf{R}^{m \times n}$ such that
\begin{equation}
f(z + d) = f(z) + \mathcal{J} d + \mathbf{O}(\|d\|^2) \,,
\end{equation}
then this function is differentiable at $z$, and the matrix 
$\mathcal{J}$ is the 
Jacobian matrix. Let us now define the Hessian matrix for a scalar 
function $f_i$. The 
Hessian 
matrix at $z$ is written $\mathcal{H}_i \in \mathbf{R}^{n \times n}$, 
where its 
entries are the second-order partial derivatives of $f_i$. For our 
vector 
function $f$, we define
\begin{equation}
\mathcal{H} = \begin{bmatrix}
\mathcal{H}_1 \\
\mathcal{H}_2 \\
\vdots \\
\mathcal{H}_m
\end{bmatrix} \in \mathbf{R}^{mn \times n} \,.
\end{equation}
This allows us to write the second-degree Taylor expansion:
\begin{equation}
f(z + d) = f(z) + \mathcal{J} d + \frac{1}{2}(I_m \otimes d^T) 
\mathcal{H} d
+ \mathbf{O}(\|d\|^3) \,,
\end{equation}
where the quadratic term is simply

\begin{equation} \label{eq:taylorsecond}
(I_m \otimes d^T) \mathcal{H} d = \begin{bmatrix} d^T \mathcal{H}_1 d 
\\ d^T 
\mathcal{H}_2 
d \\ \vdots \\ d^T \mathcal{H}_m d \end{bmatrix}^T \,.
\end{equation}

In our work with sensitivities we need to use partial 
derivatives extensively. For a function $f(x, y, p)$, where $x, y, p$ 
are vectors, we define $F_x$ to be the portion of the Jacobian matrix of
$f(x, y, p)$ w.r.t. $x$, $F_y$ the portion of the Jacobian matrix of $f(x, y, 
p)$
w.r.t. $y$, and so  on. We will often use the notation $F_{\alpha x}$, where 
$\alpha$ is an element of $p$. Because $\alpha$ is a scalar, this is equivalent 
to the matrix $\frac{dF_x}{d\alpha}$.

\subsection{Implementation of Trajectory Sensitivities}
\label{section:trajimp}
Let  us consider the DAE system (\ref{eq:dae}) and elements of the parameter 
vector p: 
$\alpha \triangleq p_i$, $\beta \triangleq p_j$. The first-order sensitivity 
vector 
with respect to $\alpha$ is written

\begin{equation}
	u^{\alpha} = \begin{bmatrix}
		u^{\alpha}_x & u^{\alpha}_y
	\end{bmatrix} =
	\begin{bmatrix}
		\frac{\partial x_1}{\partial \alpha},
		\frac{\partial x_2}{\partial \alpha},
		\dots,
		\frac{\partial y_1}{\partial \alpha},
		\dots
	\end{bmatrix} \,.
\end{equation}

The continuous first-order sensitivities can be computed as follows:

\begin{equation} \label{eq: firstsens}
\begin{bmatrix}
\dot{u}_x^{\alpha} \\
0
\end{bmatrix} = 
\begin{bmatrix}
F_x & F_y \\
G_x & G_y
\end{bmatrix}
\begin{bmatrix}
u_x^{\alpha} \\
u_y^{\alpha}
\end{bmatrix} +
\begin{bmatrix}
F_{\alpha} \\
G_{\alpha}
\end{bmatrix} \,.
\end{equation}

Equation \eqref{eq: firstsens} is a nonhomogeneous linear differential equation and can be integrated with a backward Euler scheme. If we set $\Delta t$ to be our time step, then we just need to solve the following linear system:
\begin{equation} \label{eq: firstsens_sol}
\begin{bmatrix}
\Delta t F_x - I & \Delta t F_y \\
G_x & G_y
\end{bmatrix}
\begin{bmatrix}
\{u_x^{\alpha}\}^t \\
\{u_y^{\alpha}\}^t
\end{bmatrix} =
\begin{bmatrix}
-\Delta t(F_{\alpha}) - \{u^{\alpha}\}^{t-1}x \\
-G_{\alpha}
\end{bmatrix} \,.
\end{equation}

The second-order sensitivity vector is defined as follows:

\begin{equation}
	v^{\alpha \beta} = \begin{bmatrix}
		v^{\alpha \beta}_x & v^{\alpha \beta}_y
	\end{bmatrix} =
	\begin{bmatrix}
		\frac{\partial x_1}{\partial \alpha \beta},
		\frac{\partial x_2}{\partial \alpha \beta},
		\dots,
		\frac{\partial y_1}{\partial \alpha \beta},
		\dots
	\end{bmatrix} \,.
\end{equation}

For the second-order sensitivities the structure of the equations is also a nonhomogeneous linear differential equation but has a more complex forcing term $\xi$ that depends on the first-order sensitivities:
\begin{equation} \label{eq: secondsens}
\begin{bmatrix}
\dot{v}_x^{\alpha \alpha} \\
0
\end{bmatrix} = 
\begin{bmatrix}
F_x & F_y \\
G_x & G_y
\end{bmatrix}
\begin{bmatrix}
v^{\alpha \alpha}_x \\
v^{\alpha \alpha}_y
\end{bmatrix} + \xi(x, y, u) \,,
\end{equation}
where the forcing term $\xi = \begin{bmatrix} \xi_x & \xi_y \end{bmatrix}^T$ is
\begin{equation}
\begin{bmatrix}
F_ {\alpha \alpha} \\
G_ {\alpha \alpha} 
\end{bmatrix} +
2 \begin{bmatrix}
F_{\alpha x} & F_{\alpha y}  \\
G_{\alpha x} & G_{\alpha y} 
\end{bmatrix}
\begin{bmatrix}
u^{\alpha}_x \\
u^{\alpha}_y
\end{bmatrix} + (I_m \otimes (u^{\alpha})^T) \mathcal{H} u^{\alpha} \,.
\end{equation}
Since this forcing term is given, the structure of the integration is exactly the same as before:

\begin{equation} \label{eq: secondsens_sol}
\begin{bmatrix}
\Delta tF_x - I & \Delta t(F_y) \\
G_x & G_y
\end{bmatrix}
\begin{bmatrix}
	\{v_x^{\alpha\alpha}\}^t \\
	\{v_y^{\alpha\alpha}\}^t
\end{bmatrix} =
\begin{bmatrix}
-\Delta t(\xi_x) - \{v^{\alpha \alpha}_x\}^{t-1} \\
-\xi_y
\end{bmatrix} \,.
\end{equation}

The mixed sensitivities $v^{\alpha \beta}$ are obtained with the same 
equation \eqref{eq: secondsens} but with a different 
forcing term $\xi$:

\begin{align}
	\begin{bmatrix}
		F_ {\alpha \beta} \\
		G_ {\alpha \beta} 
	\end{bmatrix} &+
	\begin{bmatrix}
		F_{\alpha x} & F_{\alpha y}  \\
		G_{\alpha x} & G_{\alpha y} 
	\end{bmatrix}
	\begin{bmatrix}
		u^{\alpha}_x \\
		u^{\alpha}_y
	\end{bmatrix} \\
	& + \begin{bmatrix}
	F_{\beta x} & F_{\beta y}  \\
	G_{\beta x} & G_{\beta y} 
\end{bmatrix}
\begin{bmatrix}
	u^{\beta}_x \\
	u^{\beta}_y
\end{bmatrix}
+ (I_m \otimes (u^{\alpha})^T) \mathcal{H} u^{\beta} \,.
\end{align}

We define  the matrix of columns $U = \begin{bmatrix}
	u^{\alpha} & u^{\beta} & \dots
\end{bmatrix}$, where each rows $u_i$ will represent the sensitivities of the 
$i$th variable with respect to all the parameters. At the same time we define 
the matrix.

\begin{equation}
	V_{x} = \begin{bmatrix}
		\frac{\partial^2 x}{\partial \alpha \alpha} & \frac{\partial^2 
		x}{\partial 
		\alpha 
		\beta} & \dots \\
		\frac{\partial^2 x}{\partial \beta \alpha} & \frac{\partial^2 
		x}{\partial 
		\beta 
		\beta} & \dots \\
		\dots & \dots & \ddots
	\end{bmatrix}
\end{equation}

For a system of $p$ parameters we will have  $p$ $n-$dimensional vectors of 
first-order sensitivities, $p$ $n-$dimensional vectors of second-order self-sensitivities, and $\frac{p^2 - p}{2}$ $n-$dimensional vectors of second-order 
mixed sensitivities. A priori, it might seem that the computations involved as 
we increase the dimension of our parameters are too onerous. However, one 
should note that the linear system that we solve for each sensitivity system is 
the same  and thus one can leverage a solver 
with multiple right-hand sides. In any case, the polynomial nature of the 
growth, compared 
with the 
exponential growth required to sample the parameter space in Monte Carlo 
methods, makes this method much more computationally feasible.

\section*{Acknowledgment}
This material was based upon work supported by the U.S. Department of Energy, Office of Science, under Contract DE-AC02-06CH11347.
\bibliographystyle{ieeetr}
\bibliography{bibliography}

\begin{thebibliography}{10}

\bibitem{Chaspierre2018}
G.~Chaspierre, P.~Panciatici, and T.~V. Cutsem, ``Modelling active distribution
  networks under uncertainty: Extracting parameter sets from randomized dynamic
  responses,'' in {\em 2018 Power Systems Computation Conference ({PSCC})},
  {IEEE}, June 2018.

\bibitem{NERC2016}
``{Reliability Guideline: Developing Load Model Composition Data},'' tech.
  rep., North American Electric Reliability Corporation, March 2017.

\bibitem{Milano2013}
F.~Milano and R.~Zarate-Minano, ``A systematic method to model power systems as
  stochastic differential algebraic equations,'' {\em {IEEE} Transactions on
  Power Systems}, vol.~28, pp.~4537--4544, Nov. 2013.

\bibitem{Dhople2013}
S.~V. Dhople, Y.~C. Chen, L.~DeVille, and A.~D. Dominguez-Garcia, ``{Analysis
  of Power System Dynamics Subject to Stochastic Power Injections},'' {\em IEEE
  Transactions on Circuits and Systems I: Regular Papers}, vol.~60,
  pp.~3341--3353, dec 2013.

\bibitem{Adeen2021}
M.~Adeen and F.~Milano, ``On the impact of auto-correlation of stochastic
  processes on the transient behavior of power systems,'' {\em {IEEE}
  Transactions on Power Systems}, vol.~36, pp.~4832--4835, Sept. 2021.

\bibitem{Roberts2016}
C.~Roberts, E.~M. Stewart, and F.~Milano, ``{Validation of the
  Ornstein--Uhlenbeck process for load modeling based on µPMU measurements},''
  in {\em 2016 Power Systems Computation Conference (PSCC)}, pp.~1--7, IEEE,
  jun 2016.

\bibitem{Dong2012}
Z.~Y. Dong, J.~H. Zhao, and D.~J. Hill, ``{Numerical Simulation for Stochastic
  Transient Stability Assessment},'' {\em IEEE Transactions on Power Systems},
  vol.~27, pp.~1741--1749, nov 2012.

\bibitem{Hockenberry2004}
J.~Hockenberry and B.~Lesieutre, ``Evaluation of uncertainty in dynamic
  simulations of power system models: The probabilistic collocation method,''
  {\em {IEEE} Transactions on Power Systems}, vol.~19, pp.~1483--1491, Aug.
  2004.

\bibitem{Xu2019}
Y.~Xu, L.~Mili, A.~Sandu, M.~R. von Spakovsky, and J.~Zhao, ``Propagating
  uncertainty in power system dynamic simulations using polynomial chaos,''
  {\em {IEEE} Transactions on Power Systems}, vol.~34, pp.~338--348, Jan. 2019.

\bibitem{jovica}
J.~V. Milanovi{\'{c}}, ``Probabilistic stability analysis: the way forward for
  stability analysis of sustainable power systems,'' {\em Philosophical
  Transactions of the Royal Society A: Mathematical, Physical and Engineering
  Sciences}, vol.~375, p.~20160296, July 2017.

\bibitem{osti_1644285}
A.~Faris, D.~Kosterev, J.~H. Eto, and D.~Chassin, ``Load composition analysis
  in support of the {NERC Load Modeling Task Force} 2019--2020 field test of
  the composite load model,'' 6 2020.

\bibitem{Hiskens2006}
I.~A. Hiskens and J.~Alseddiqui, ``Sensitivity, approximation, and uncertainty
  in power system dynamic simulation,'' {\em IEEE Transactions on Power
  Systems}, vol.~21, no.~4, pp.~1808--1820, 2006.

\bibitem{Kim2020}
J.-K. Kim, B.~Lee, J.~Ma, G.~Verbic, S.~Nam, and K.~Hur, ``Understanding and
  evaluating systemwide impacts of uncertain parameters in the dynamic load
  model on short-term voltage stability,'' {\em {IEEE} Transactions on Power
  Systems}, pp.~1--1, 2020.

\bibitem{Schweppe_Fred_C1973}
F.~C. Schweppe, {\em Uncertain dynamic systems}.
\newblock Prentice-Hall, hardcover~ed., 1973.

\bibitem{lesieutrereport}
B.~Lesieutre, ``{Improving Dynamic Load and Generator Response Performance
  Tools},'' tech. rep., LBNL, Berkeley, 2005.

\bibitem{nerccriteria}
NERC, ``{Standard TPL-001-4 — Transmission System Planning Performance
  Requirements}.'' \url{https://www.nerc.com/files/TPL-001-4.pdf}.

\bibitem{tenza2016analysis}
M.~Tenza and S.~Ghiocel, ``An analysis of the sensitivity of {WECC} grid
  planning models to assumptions regarding the composition of loads,'' {\em
  Mitsubishi Electric Power Products}, 2016.

\bibitem{Timko1983}
K.~Timko, A.~Bose, and P.~Anderson, ``{Monte Carlo} simulation of power system
  stability,'' {\em {IEEE} Transactions on Power Apparatus and Systems},
  vol.~{PAS}-102, pp.~3453--3459, Oct. 1983.

\bibitem{Zhang2020}
Y.~Zhang, Y.~Li, K.~Tomsovic, S.~M. Djouadi, and M.~Yue, ``Review on
  set-theoretic methods for safety verification and control of power system,''
  {\em {IET} Energy Systems Integration}, vol.~2, pp.~226--234, Sept. 2020.

\bibitem{Zhang2017}
H.~Zhang, S.~Abhyankar, E.~Constantinescu, and M.~Anitescu, ``Discrete adjoint
  sensitivity analysis of hybrid dynamical systems with switching,'' {\em IEEE
  Transactions on Circuits and Systems I: Regular Papers}, vol.~64, no.~5,
  pp.~1247--1259, 2017.

\bibitem{Geng2019}
S.~Geng and I.~A. Hiskens, ``Second-order trajectory sensitivity analysis of
  hybrid systems,'' {\em IEEE Transactions on Circuits and Systems I: Regular
  Papers}, vol.~66, no.~5, pp.~1922--1934, 2019.

\bibitem{Choi2017}
H.~Choi, P.~J. Seiler, and S.~V. Dhople, ``Propagating uncertainty in
  power-system {DAE} models with semidefinite programming,'' {\em IEEE
  Transactions on Power Systems}, vol.~32, no.~4, pp.~3146--3156, 2017.

\bibitem{Milano2010}
F.~Milano, {\em Power System Modelling and Scripting}.
\newblock Springer Berlin Heidelberg, 2010.

\bibitem{Hiskens2000}
I.~Hiskens and M.~Pai, ``{Trajectory sensitivity analysis of hybrid systems},''
  {\em IEEE Transactions on Circuits and Systems I: Fundamental Theory and
  Applications}, vol.~47, no.~2, pp.~204--220, 2000.

\bibitem{Nocedal2006Numerical}
J.~Nocedal and S.~J. Wright, {\em {Numerical Optimization}}.
\newblock Springer Series in Operations Research and Financial Engineering, New
  York: Springer, 2006.

\bibitem{Conn2000}
A.~R. Conn, N.~I.~M. Gould, and P.~L. Toint, {\em Trust Region Methods}.
\newblock Society for Industrial and Applied Mathematics, Jan. 2000.

\bibitem{lindon_roberts}
``\textit{Trustregion}: Trust-region subproblem solver.''
  \url{https://github.com/lindonroberts/trust-region}.

\bibitem{Gould1999}
N.~I.~M. Gould, S.~Lucidi, M.~Roma, and P.~L. Toint, ``Solving the trust-region
  subproblem using the {L}anczos method,'' {\em {SIAM} Journal on
  Optimization}, vol.~9, pp.~504--525, Jan. 1999.

\end{thebibliography}

\end{document}